\documentclass[12pt,a4paper]{article}
\usepackage{fullpage}
\usepackage{latexsym}
\usepackage{amssymb}
\DeclareFontEncoding{OT2}{}{}
\DeclareFontSubstitution{OT2}{wncy}{m}{n}
\DeclareFontFamily{OT2}{wncy}{}
\DeclareFontShape{OT2}{wncy}{m}{n}{<5-9> gen * wncyr <10-> wncyr10}{}
\DeclareFontShape{OT2}{wncy}{m}{sc}{<-> wncysc10}{}
\DeclareSymbolFont{wncyr}{OT2}{wncy}{m}{n}
\DeclareMathSymbol\DL{\mathord}{wncyr}{68}

\def\D{\Delta}
\def\Q{\Delta\hbox to 0pt{\hss$\scriptstyle\Delta$}}
\def\QQ{\Delta\hbox to 0pt{\hss$\scriptscriptstyle\Delta$}}

\def\R{\mathbb{R}}
\def\S{\mathop{\mathrm{Sym}}}
\def\#{\sharp}

\def\j{\mathrm{jet}}
\def\J{\mathop{\mathrm{Jet}}}

\def\ev{\mathrm{ev}\,}
\def\E{\mathcal{E}}
\def\End{\mathrm{End}\,}
\def\SO{\mathbf{SO}\,}
\def\[#1]{[\![\,#1\,]\!]}
\def\<#1,#2>{\langle\,#1,\,#2\,\rangle}
\def\pskip{\par\vskip2mm plus1mm minus1mm}
\newtheorem{Proposition}{Proposition}[section]
\newtheorem{Lemma}[Proposition]{Lemma}
\newtheorem{Theorem}[Proposition]{Theorem}

\newtheorem{Definition}[Proposition]{Definition}
\newtheorem{Example}[Proposition]{Example}
\newtheorem{Remark}[Proposition]{Remark}
\newenvironment{Reference}[1]{\par\pskip\noindent\textbf{#1}\bgroup\em}
 {\egroup\par\pskip}
\def\endsummary{\hfil\break}
\newenvironment{Proof}{\par\noindent\textbf{Proof:\ }}{\hfill
 \ensuremath{\Box}\par}
\begin{document}
\title{A Characterization of the Heat Kernel Coefficients}
\author{Gregor Weingart\footnote{partially supported by the SFB 256
 \textsl{``Nichtlineare partielle Differentialgleichungen''}}}
\def\today{May 16, 2001}
\maketitle
\hyphenation{ma-ni-fold}
\begin{abstract}
 We consider the asymptotic expansion of the heat kernel of
 a generalized Laplacian for $t\rightarrow 0^+$ and characterize
 the coefficients $a_k,\,k\geq 0,$ of this expansion by a natural
 intertwining property. In particular we will give a closed formula
 for the infinite order jet of these coefficients on the diagonal in
 terms of the local expressions of the powers of the given generalized
 Laplacian in normal coordinates.
\end{abstract}
\begin{center}
\end{center}
\section{Introduction}
 A generalized Laplacian $\Q$ acting on sections of a vector bundle $\E$
 over a Riemannian manifold $M$ is a second order elliptic differential
 operator with scalar symbol (\cite{bgv}). Its spectrum is interesting
 for various geometrical questions and the spectra of the more prominent
 generalized Laplacians are of fundamental importance in global analysis.
 The short distance behaviour of the spectrum of a generalized Laplacian
 can be studied using analytic properties of its heat operator $K^{\QQ}_t$,
 which solves the heat equation
 $$
  \frac\partial{\partial t}\,(\,K^{\QQ}_t\psi\,)\quad=\quad-\;\Q\psi
 $$
 associated to $\Q$ for all sections $\psi$ of $\E$. The general existence
 result of Minakshisundaram--Pleijel (\cite{gil},\cite{bgv}) shows that
 $K^{\QQ}_t$ is an integral operator for all $t>0$ with integral kernel
 $k^{\QQ}_t(x,y)$, which is a smooth section of the bundle with fiber
 $\mathrm{Hom}(\E_x,\E_y)$ over the point $(x,y)$ of the product $M\times M$.
 Except for a very limited set of examples the heat kernel is not known
 explicitly, perhaps the most important exception is the heat kernel of
 the flat Laplacian $\D$ acting on sections of a trivial vector bundle
 $V\times E$ over a euclidian vector space $V$ of dimension $n$:
 $$
  k^{\D}_t(x,y)\quad=\quad
  \frac1{\sqrt{4\pi t}^n}\,e^{\textstyle-\frac{|x-y|^2}{4t}}\,\mathrm{id}_E
 $$
 As can be seen from this example the heat kernel $k^{\QQ}_t(x,y)$ peaks
 on the diagonal $\{x=y\}$ of $M\times M$ as $t\rightarrow 0^+$ and
 it turns out that the rate of its divergence on the diagonal carries a
 tremendous amount of information about the geometry of the underlying
 manifold $M$. More precisely the heat kernel has an asymptotic expansion
 along the diagonal as $t\rightarrow 0^+$ of the form
 \begin{equation}\label{ae}
  k^{\QQ}_t(x,y)\quad\vtop{\hbox{$\displaystyle\sim$}
   \vskip-2mm\hbox{$\scriptstyle t\rightarrow 0^+$}}\quad
  \frac1{\sqrt{4\pi t}^n}\,e^{\textstyle-\frac{\mathrm{dist}^2(x,y)}{4t}}
  \,\sum_{k\geq 0}t^\mu\,\hat a_k(x,y)
 \end{equation}
 where $\mathrm{dist}^2(x,y)$ is the square of the geodesic distance from
 $x$ to $y$ and the $\hat a_k(x,y)$ are again sections of the bundle
 $\mathrm{Hom}(\E_x,\E_y)$. In a sense these coefficients $\hat a_k,
 \,k\geq 0,$ describe how much the heat kernel of the operator $\Q$
 differs in divergence from the flat Laplacian $\D$ in various directions
 and this difference does not only depend on the coefficients of $\Q$, but
 involves the geometry of $M$ in a crucial way.

 As a matter of convenience we will not use the coefficients $\hat a_k(x,y)$
 directly but will multiply them $a_k(x,y)\,:=\,j(x,y)\hat a_k(x,y)$ by the
 jacobian determinant $j(x,y)$ of the differential of the exponential map
 in $y$ at $x$, which is a well defined smooth function in a neighborhood
 of the diagonal. Off this neighborhood we may extend $j(x,y)$ arbitrarily to
 a smooth function on $M\times M$, because only the germs of the coefficients
 $\hat a_k(x,y)$ along the diagonal are characterized by the asymptotic
 expansion (\ref{ae}) of $k^{\QQ}_t$ above. The main theorem of this article
 exhibits a remarkable intertwining property of these modified coefficients
 $a_k,\,k\geq 0$:
 
 \begin{Reference}{Theorem \ref{iprop}}
 {\rm\textbf{(The Intertwining Property)}}\endsummary
  Let $\Q$ be a generalized Laplacian acting on sections $\psi$ of a vector
  bundle $\E$ over $M$. The coefficients $a_k:\E_x\longrightarrow\E_y$ of
  the asymptotic expansion of the heat kernel of $\Q$ satisfy
  $$
   \frac{(-1)^\mu}{\mu!}\,\left[\,\Q^\mu\,\psi\,\right](y)
   \quad=\quad
   \sum_{\nu=0}^\mu\,\frac{(-1)^\nu}{\nu!}\,\left[\,\D^\nu\,(a_{\mu-\nu}
   (\mathrm{exp}_y\cdot,y)\,\psi(\mathrm{exp}_y\cdot)\,\right](0)
  $$
  where $a_{\mu-\nu}(\mathrm{exp}_y\cdot,y)\,\psi(\mathrm{exp}_y\cdot)$ is
  considered as a section of the trivial vector bunde $T_yM\times\E_y$ over
  $T_yM$ and $\D$ is the flat Laplacian acting on sections of this trivial
  bundle.
 \end{Reference}

 The rest of the article is devoted to corollaries of Theorem \ref{iprop}. In
 particular we prove the striking fact that the intertwining property alone
 characterizes the infinite order jet of all coefficients $a_k(\cdot,y)$ in
 $y$ and on the way we will give a different proof of Polterovich's inversion
 formula (\cite{polt1}) for the endomorphisms $a_k(y,y),\,k\geq 0,$ of the
 fiber $\E_y$ in a somewhat stronger formulation:

 \begin{Reference}{Theorem \ref{ip}}
 {\rm\textbf{(Polterovich's Inversion Formula)}}\endsummary
  Let $\Q$ be a generalized Laplacian acting on sections of a vector bundle
  $\E$ over $M$ and let $\mathrm{dist}^2(\cdot,y)$ be the square of the
  geodesic distance to $y\in M$. For all $r\geq k\geq 0$ and all sections
  $\psi$ of $\E$ the endomorphism $a_k(y,y)\,=\,\hat a_k(y,y)$ of the
  fiber $\E_y$ is given by the formula:
  $$
  \left[\,a_k\,\psi\,\right](y)\quad=\quad
  \sum_{l=0}^r\,\left(-{\textstyle\frac14}\right)^l\,{r+\frac{n}2\choose r-l}
  \,\left[\,\frac{(-1)^{k+l}}{(k+l)!}\,\Q^{k+l}\,
  (\,\frac1{l!}\,\mathrm{dist}^{2l}(\cdot,y)\,\psi\,)\,\right](y)\,.
  $$
 \end{Reference}
 
 The final result is Theorem \ref{dlf}, which provides an explicit formula
 for the infinite order jets of the coefficients $a_k(\cdot,y)$ in $y$ in
 terms of the local expressions of the powers of the operator $\Q$ in normal
 coordinates on $M$ and some smooth trivialization of $\E$.

 \pskip
 The present version of this article was written down while I was staying at
 the University of California in Berkeley, whose hospitality is gratefully
 acknowledged. It is only the first part of a projected article of three and
 hence will not appear separately in printed form.
\section{The Intertwining Property}
 It is well known that the coefficients $\hat a_k,k\geq 0,$ of the asymptotic
 expansion of the heat kernel of $\Q$ are local expressions in the curvature
 of $M$ and the coefficients of $\Q$ (\cite{gil}). Consequently we will assume
 throughout this section that $\Q$ is defined on sections of a trivializable
 vector bundle $\E$ over a starshaped neighborhood of the origin in some
 euclidean vector space $V$ of dimension $n$. Moreover we will assume that
 the identity map $V\longrightarrow V$ is the exponential map in the origin
 with respect to the Riemannian metric defined by the symbol of $\Q$, a 
 fortiori the symbol metric must agree with the scalar product of $V$ in $0$.
 The general case is easily reduced to this local model by taking $V\,=\,T_yM$
 to be the tangent space of the underlying manifold at the point $y$ under
 consideration and choosing the neighborhood of the origin and the
 trivializable vector bundle $\E$ accordingly.

 The preceeding assumptions are crucial to the entire argument but single
 out the origin as a special point. Consequently all our statements below
 refer to the fixed target point $y=0$. In particular the coefficients
 $\hat a_k(x):=\hat a_k(x,0)$ will be homomorphisms $\hat a_k(x):\;\E_x
 \longrightarrow E$ from the fiber of $\E$ over $x$ to the fiber $E$ of
 $\E$ over the origin. Similarly we will write $j(x):=j(x,0)$ for the
 Jacobian determinant of the exponential map in the origin at the point
 $x\in V$, which is the identity map $V\longrightarrow V$ by assumption.
 In other words the Jacobian determinant is characterized by $d\mathrm
 {vol}_g\,=\,j(x)dx$ and is the square root $j(x)\,:=\,\sqrt{\mathrm{det}
 \,g(x)}$ of the determinant of the symbol metric with respect to the flat
 scalar product. Now by the remarks following (\cite[Theorem 2.30]{bgv})
 we have an asymptotic expansion for the heat operator
 $$
  \left[\,K^{\QQ}_t\psi\,\right](0)\quad\vtop{\hbox{$\displaystyle\sim$}
   \vskip-2mm\hbox{$\scriptstyle t\rightarrow 0^+$}}\quad
  \sum_{\mu\geq 0}\;\frac{(-t)^\mu}{\mu!}\,\left[\,\Q^\mu\psi\,\right](0)
 $$
 which complements the asymptotic expansion (\ref{ae}) of its integral
 kernel $k^{\QQ}_t$. Using both asymptotic expansions we find for any
 smooth section $\psi$ of $\E$ with compact support in a sufficiently
 small neighborhood of the origin
 \begin{eqnarray*}
 \int_V\,k^{\QQ}_t(x)\psi(x)\,d\,\mathrm{vol}_g
  &\vtop{\hbox{$\displaystyle\sim$}\vskip-2mm
    \hbox{$\scriptstyle t\rightarrow 0^+$}}&
   \sum_{\mu\geq 0}\:\frac{(-t)^\mu}{\mu!}\,\left[\,\Q^\mu\psi\,\right](0)\\
  &\vtop{\hbox{$\displaystyle\sim$}\vskip-2mm
    \hbox{$\scriptstyle t\rightarrow 0^+$}}&
   \int_V\,\frac1{\sqrt{4\pi t}^n}\,e^{\textstyle-\frac{|x|^2}{4t}}\,
    \sum_{\mu\geq 0}\,t^\mu\,\hat a_\mu(x)\,\psi(x)\,d\mathrm{vol}_g\\
  &\vtop{\hbox{$\displaystyle\sim$}\vskip-2mm
    \hbox{$\scriptstyle t\rightarrow 0^+$}}&
   \sum_{\mu\geq 0}\,t^\mu\,\int_V\frac1{\sqrt{4\pi t}^n}\,
    e^{\textstyle-\frac{|x|^2}{4t}}\,(\,j(x)\,\hat a_\mu(x)\,\psi(x)\,)\,dx\\
  &\vtop{\hbox{$\displaystyle\sim$}\vskip-2mm
    \hbox{$\scriptstyle t\rightarrow 0^+$}}&
   \sum_{\mu\geq 0}\,t^\mu\,\sum_{\nu\geq 0}\,\frac{(-t)^\nu}{\nu!}
    \,\left[\,\D^\nu(\,a(x)\,\psi(x)\,)\,\right](0)\,,
 \end{eqnarray*}
 where $a(x)\,:=\,j(x)\hat a_k(x)$ and $\D$ is the flat Laplacian acting
 on sections of the trivial $E$--bundle $V\times E$ over $V$. Note that the
 asymptotic expansion of the heat kernel $k^{\QQ}_t$ is locally uniform in
 $x$ according to (\cite[Theorem 2.30]{bgv}) and hence we may integrate it
 over the compact support of $\psi$ to obtain an asymptotic expansion of
 the integral. Existence of asymptotic expansions implies uniqueness and so
 we need only sort out the different powers of $t$ to prove the intertwining
 property of the heat kernel coefficients in the form:
 
 \begin{Theorem}[Intertwining Property of the Heat Kernel Coefficients]
 \label{iprop}\endsummary
  Let $\Q$ be a generalized Laplacian acting on sections of a trivializable
  vector bundle $\E$ over $V$. The coefficients $a_k,\,k\geq 0,$ of the
  asymptotic expansion of the heat kernel of $\Q$ are sections of the bundle
  $\mathrm{Hom}(\E,\,V\times E)$ and so $a_k\psi$ is a well defined section
  of $V\times E$ for any section $\psi$ of $\E$. The coefficients $a_k$
  intertwine the powers of $\Q$ and $\D$:
  \begin{equation}\label{inter}
   \frac{(-1)^\mu}{\mu!}\,\left[\,\Q^\mu\,\psi\,\right](0)
   \quad=\quad
   \sum_{\nu=0}^\mu\,\frac{(-1)^\nu}{\nu!}\,\left[\,\D^\nu\,(a_{\mu-\nu}\,
   \psi)\,\right](0)\,.
  \end{equation}
 \end{Theorem}

 A direct consequence of Theorem \ref{iprop} is that a generalized Laplacian
 $\Q$ acting on sections of a vector bundle $\E$ remembers to vanish on
 sections $\psi$ of $\E$, which are harmonic in some smooth trivialization
 of $\E$, in the following sense:

 \begin{Remark}[Excess Vanishing for Harmonic Sections]
 \endsummary
  Consider a section $\psi$ of $\E$, which is a homogeneous, harmonic
  polynomial $\Phi\psi\,\in\,\S^dV^*\otimes E$ of degree $d$ with $\D(\Phi\psi)
  \,=\,0$ in some smooth trivialization $\Phi$ of $\E$. If $\Q$ is a
  generalized Laplacian acting on sections of $\E$, then $\left[\,
  \Q^\mu\psi\,\right](0)\,=\,0$ vanishes for all $\mu<d$.
 \end{Remark}

 Philosophically the intertwining property (\ref{inter}) provides a precise
 geometric interpretation for the coefficients $a_k$ of the heat kernel
 expansion in terms of $\Q$ and the flat model operator $\D$ for generalized
 Laplacians. Similar considerations should apply to other model operators
 arising e.~g.~in Heisenberg calculus or in other parabolic calculi. A
 convenient reformulation of this property (\ref{inter}) can be given
 with the help of the formal power series $e^{-z\Q}$ and $e^{-z\D}$ of
 differential operators on $\E$ and $V\times E$ respectively and the
 generating series
 $$
  a(z)\quad:=\quad \sum_{\mu\geq 0}\;a_\mu\,z^\mu
 $$
 for the heat kernel coefficients $a_k,\,k\geq0$. Namely the intertwining
 property (\ref{inter}) is just another way to write down the equality
 $[e^{-z\Q}\psi](0)\,=\,[e^{-z\D}(a(z)\psi)](0)$ of formal power series
 in $z$ for all sections $\psi$, which in turn is equivalent to the
 commutativity of the pentagram
 \begin{center}\unitlength.5mm
  \begin{picture}(166,80)
   \put(  8, 67){$\J_0^\infty\E$}
   \put( 34, 70){\vector(1,0){89}}
   \put(125, 67){$\J_0^\infty E\[z]$}
   \put( 19, 63){\vector(0,-1){28}}
   \put(  4, 27){$\J_0^\infty \E\[z]$}
   \put(144, 63){\vector(0,-1){28}} 
   \put(125, 27){$\J_0^\infty E\[z]$}
   \put( 20, 24){\vector(3,-1){47}}
   \put(137, 24){\vector(-3,-1){47}}
   \put( 68,  3){$E\[z]$}
   \put( 71, 73){$a(z)$}
   \put(  0, 47){$e^{-z\QQ}$}
   \put(147, 47){$e^{-z\D}$}
   \put( 34, 10){$\ev$}
   \put(114, 10){$\ev$}
  \end{picture}
 \end{center}
 where $\ev$ is the evalutation at the origin. Only the infinite order jet
 of the sections $a_k,\,k\geq 0,$ of $\Gamma(V,\mathrm{Hom}(\E,V\times E))$
 in the origin can ever be sensed by jet theory of course, and so the
 commutativity of this square does only depend on the infinite order
 jet of the generating series $a(z)$. Strikingly however the commutativity
 of the pentagon characterizes the infinite order jet of $a(z)$ completely.
 In fact the unique jet solution of the intertwining property thought of as
 a set of equations for the unknowns $a_k,\,k\geq 0,$ is given in Theorem
 \ref{dlf} below.
 
 \pskip
 At this point we digress a little bit on a very interesting property of
 the flat Laplacian $\D$ acting on functions $C^\infty V$ on a euclidian
 vector space $V$. More general we can consider the flat Laplacian $\D$
 acting on sections $C^\infty(V,E)$ of the trivial $E$--bundle $V\times E$
 over $V$, but the auxiliary vector space $E$ never enters into the argument
 directly and so we will stick to the case $E=\R$. Along with the flat
 Laplacian $\D$ comes the operator $|x|^2$ of multiplication by the
 square of the distance to the origin, in orthogonal coordinates
 $\{x_\mu\}$ on $V$ these two operators can be written:
 $$
  \D\quad:=\quad-\,\sum_\mu\,\frac{\partial^2}{\partial x_\mu^2}
  \qquad\qquad\qquad
  |x|^2\quad:=\quad\sum_\mu\,x_\mu^2\,.
 $$
 By restriction they act on the subspace $\S\,V^*\subset C^\infty V$ of
 polynomials on $V$ and their commutator $[\D,\,|x|^2]\,=\,(-4)\,(N+\frac{n}2)$
 as operators on $\S\,V^*$ involves the number
 operator $N:\;\S\,V^*\longrightarrow\S\,V^*$, which multiplies $\S^rV^*$
 by $r$, shifted by half the dimension $n$ of $V$. In other words
 $X:=\frac12|x|^2$ and $Y:=\frac12\D$ close with $H:=N+\frac{n}2$
 to an $\mathfrak{sl}_2$--algebra of operators on $\S\,V^*$. Iterated
 commutators of $X$ and $Y$ in $\mathfrak{sl}_2$--representations can
 in general be written down using factorial polynomials $[x]_r\,:=\,
 x(x-1)\ldots(x-r+1)$ or binomial coefficients ${x\choose r}\,:=\,
 \frac1{r!}[x]_r$ for integral $r\geq 0$. In particular the standard
 relation $Y^r\,X^r\,\psi\,=\,r!\,[-\lambda]_r\,\psi$ for all $\psi$
 with $Y\,\psi\,=\,0$ and $H\,\psi\,=\,\lambda\,\psi$ implies for the
 constant function $\psi\,:=\,1\,\in\,\S^0V^*$ the classical formula:
 \begin{equation}\label{clp}
  \left[\,\D^r\,|x|^{2r}\,\right](0)\quad=\quad 4^r\;r!\;
  [\,-{\textstyle\frac{n}2}\,]_r\,.
 \end{equation}
 Slightly more useful for our purposes is the following derived identity:

 \begin{Lemma}\label{mi}
  For all smooth functions $\psi\in C^\infty V$ and all $k,\,l\,\geq\,0$:
  $$
   \left[\,\frac{(-1)^{k+l}}{(k+l)!}\,\D^{k+l}\,
   (\,\frac1{k!}\,|x|^{2k}\,\psi\,)\,\right](0)
   \quad=\quad
   (-4)^k\,{-\frac{n}2-l\choose k}\,\left[\,\frac{(-1)^l}{l!}\D^l\,
   \psi\,\right](0)\,.
  $$
 \end{Lemma}

 \begin{Proof}
 Only a finite number of partial derivates of $\psi$ in the origin $0$ are
 actually involved in this identity and hence we may assume that $\psi$
 is a polynomial without loss of generality. Moreover only the homogeneous
 component of $\psi$ of degree $2l$ contributes to left and right hand
 side, which are both evidently $\SO V$--invariant linear functionals in
 $\psi\,\in\,\S^{2l}V^*$. However there is but one $\SO V$--invariant
 linear functional on $\S^{2l}V^*$ up to scale, so that it is sufficient
 to check the identity in question, which can be rewritten as
 $$
  \left[\,\D^{k+l}\,(\,|x|^{2k}\,\psi\,)\,\right](0)
  \quad=\quad
  \frac{4^{k+l}\,(k+l)!\,[-\frac{n}2]_{k+l}}
  {4^l\;\;\;l!\;\;\;[-\frac{n}2]_l}
  \left[\,\D^l\,\psi\,\right](0)\,,
 $$
 for the single polynomial $\psi\,:=\,|x|^{2l}$, for which it is true by
 the classical formula (\ref{clp}).
 \end{Proof}

 \pskip
 Returning to the general case of a Laplacian $\Q$ acting on sections of a
 vector bundle $\E$ over $V$ we recall that the generating series $a(z)$ for
 the coefficients in the asymptotic expansion of the heat kernel intertwines
 the formal power series $e^{-z\QQ}$ and $e^{-z\D}$ of differential operators.
 Using this intertwining property (\ref{inter}) together with Lemma \ref{mi}
 we calculate:
 \begin{eqnarray*}
  \left[\,\frac{(-1)^{k+l}}{(k+l)!}\,\Q^{k+l}\,
  (\,\frac1{l!}\,|x|^{2l}\,\psi\,)\,\right](0)
  &=&
  \sum_{\mu=0}^{k+l}\,
  \left[\,\frac{(-1)^\mu}{\mu!}\,\D^\mu\,
  (\,\frac1{l!}\,|x|^{2l}\,a_{k+l-\mu}\,\psi\,)\,\right](0)\\
  &=&
  \sum_{\mu=0}^k\,
  \left[\,\frac{(-1)^{l+\mu}}{(l+\mu)!}\,\D^{l+\mu}\,
  (\,\frac1{l!}\,|x|^{2l}\,a_{k-\mu}\,\psi\,)\,\right](0)\\
  &=&
  (-4)^l\,\sum_{\mu=0}^k\,{-\frac{n}2-\mu\choose l}\,
   \left[\,\frac{(-1)^\mu}{\mu!}\D^\mu\,(\,a_{k-\mu}\,\psi\,)
   \,\right](0)\,.
 \end{eqnarray*}
 Combining this equation with the binomial inversion formula
 $$
  \sum_{l=0}^r\,{r+\frac{n}2\choose r-l}\,
  {-\frac{n}2-\mu\choose l}\quad=\quad{r-\mu\choose r}
  \quad=\quad\delta_{\mu,\,0}
 $$
 valid as soon as $r\geq\mu\geq 0$ we eventually arrive at the following
 inversion formula for the heat kernel coefficients $a_k,\,k\geq 0,$
 of a generalized Laplacian $\Q$:

 \begin{Theorem}[Polterovich's Inversion Formula (\cite{polt1},\cite{polt2})]
 \label{ip}\endsummary
  Let $\Q$ be a generalized Laplacian acting on sections of a vector bundle
  $\E$ over $V$ and let $|x|$ be the radial distance from the origin. For any
  section $\psi\,\in\,\Gamma(V,\E)$ and all $r\geq k\geq 0$ we can compute
  the action of the endomorphism $a_k(0)$ of the fiber $E$ of $\E$ over the
  origin on $\psi$ by means of an explicit inversion formula:
  $$
  \left[\,a_k\,\psi\,\right](0)\quad=\quad
  \sum_{l=0}^r\,\left(-{\textstyle\frac14}\right)^l\,{r+\frac{n}2\choose r-l}
  \,\left[\,\frac{(-1)^{k+l}}{(k+l)!}\,\Q^{k+l}\,
  (\,\frac1{l!}\,|x|^{2l}\,\psi\,)\,\right](0)\,.
  $$
 \end{Theorem}

 Our normalization of the coefficients $a_k,\,k\geq 0,$ drops the factor
 $(4\pi)^{-\frac{n}2}$ arising from the value of the euclidian heat kernel
 for the flat Laplacian $\D$ on $V$ at the origin to have the intertwining
 property (\ref{inter}) in as simple a form as possible. In the original
 formulation of Polterovich this factor is part of $a_k$ and of course
 in every conceivable application this factor has to be reinserted by hand.

 \begin{Example}[The Local Index Theorem for the Gau\ss--Bonnet Operator]
 \label{gbc}\endsummary
  Choose a smooth Riemannian metric $g$ on $V$ satisfying our standing
  assumptions and a vector bundle $\E$ over $V$ with a connection $\nabla$.
  For a generalized Laplacian on $\E$ of the form $\Q_\hbar\,=\,\nabla^*
  \nabla\,+\,\hbar F$ for some smooth endomorphism $F$ of $\E$ we
  evidently have
  \begin{eqnarray*} 
   \frac{(-1)^{k+l}}{(k+l)!}\,(\,\nabla^*\nabla\,+\,\hbar\,F\,)^{k+l}
   &=& \frac{(-\hbar)^k}{k!}F^k\,\frac{(-1)^l}{l!}\D^l\;+\;
   \textrm{terms of lower order in $\hbar$}\\
   && \qquad\qquad\qquad\qquad
   \;+\;\textrm{differential operators of lower order}
  \end{eqnarray*}
  for all $k,\,l\geq 0$. Plugging this into Polterovich's Inversion
  formula we verify immediately
  $$
   a_k(0)\quad=\quad \frac{(-\hbar)^k}{k!}\,F^k\;+\;
   \textrm{terms of lower order in $\hbar$}
  $$
  for all $k\geq 0$, which may have been anticipated from Theorem \ref{iprop}.
  Remarkably this simple calculation together with a few general considerations
  concerning the order of the coefficients $a_k$ is already sufficient to prove
  the local index theorem for the Gau\ss--Bonnet operator!
 \end{Example}

 \pskip
 In conclusion the intertwining property (\ref{inter}) alone is sufficient
 to determine the value of the coefficients $a_k,\,k\geq 0,$ of the heat
 kernel expansion at $0$. Before we proceed to show that in fact the jets
 of infinite order $\j^\infty_0a_k,\,k\geq 0,$ are determined by (\ref{inter})
 we want to make a few general remarks concerning our guiding philosphy in
 the calculations to come. We want to avoid formulas involving compositions
 of differential operators, because such formulas are hardly if ever of any
 use in explicit calculations. In favourable situations it may still be
 possible to calculate the values $\ev\mathcal{D}^r,\,r\geq 0,$ of powers
 of a differential operator $\mathcal{D}$ in the origin without knowing
 the partial derivatives of their coefficients. Note that the usual
 arguments of symbol calculus become meaningless if we have no control
 over the partial derivatives of the differential operators in a singular
 point like the origin.

 In general the value of a differential operator $\mathcal{D}$ acting
 on sections of $\E$ at the origin will be an element of $\mathrm{Hom}
 (\,\J^\infty_0\E,E)$. For the trivial $E$--bundle $V\times E$ however
 we may identify $\J^\infty_0(V\times E)$ with the formal power series
 completion of $\S\,V^*\otimes E$ and hence $\mathrm{Hom}(\,\J^\infty_0
 (V\times E),E\,)$ with $\S\,V\otimes\End E$ in the usual way. The scalar
 product of the euclidian vector space $V$ extends to a scalar product
 on $\S\,V^*$ defined via Gram's permanent and characterized by
 $\<e^\alpha,e^\beta>_{\S\,V^*}\,=\,e^{\<\alpha,\beta>}$ for all $\alpha,
 \,\beta\,\in\,V^*$ with a slight abuse of notation. Alternatively we
 may choose orthonormal coordinates $\{x_\mu\}$ on $V$ and write down
 the scalar product directly
 \begin{eqnarray*}
  \<\psi,\tilde\psi>_{\S\,V^*}
  &:=&
  \sum_{r\geq 0}\,\frac1{r!}\,\sum_{\mu_1,\ldots,\mu_r}\,
  \left[\,\frac{\partial^r}{\partial x_{\mu_1}\ldots
  \partial x_{\mu_r}}\psi\,\right](0)\,
  \left[\,\frac{\partial^r}{\partial x_{\mu_1}\ldots
  \partial x_{\mu_r}}\tilde\psi\,\right](0)\\
  &=&
  \left(\,(\ev\otimes\ev)\,\circ\,e^{\<\nabla,\nabla>}\,\right)
  \,(\,\psi\otimes\tilde\psi\,)
 \end{eqnarray*}
 where $\<\nabla,\nabla>$ is the bidifferential operator $\psi\otimes
 \tilde\psi\longmapsto\sum_\mu\,\frac\partial{\partial x_\mu}\psi\otimes
 \frac\partial{\partial x_\mu}\tilde\psi$ and $\ev$ is the evaluation at
 $0$ as before. Written in this form it is clear that the scalar product
 extends to the formal power series completion of $\S\,V^*$ or even to smooth
 functions provided that the defining sum converges. Moreover the musical
 isomorphism $\#:\;\S\,V\longrightarrow\S\,V^*$ of this scalar product is
 the natural extension of the musical isomorphism of $V$ and extends to
 the formal power series completion of $\S\,V$, too. The value of a
 differential operator $\mathcal{D}$ acting on sections of $V\times E$
 has an image $(\ev\,\mathcal{D})^\#$ under this musical isomorphism,
 which is the unique polynomial on $V$ satisfying $\left[\mathcal{D}\psi
 \right](0)\,=\,\<(\ev\,\mathcal{D}^\#),\psi>_{\S\,V^*}$ for all smooth
 sections $\psi$ of $V\times E$.

 The scalar product on $\S\,V^*$ will enter the formulas through the
 operator $\<\!\nabla,\!\nabla\!>$ introduced above, which in turn makes
 its appearance via Green's identity for the flat Laplacian $\D$ and the
 multiplication map $m$
 $$
  (\,\D\,\circ\,m\,)\,(\,\psi\otimes\tilde\psi\,)
  \quad=\quad(\,m\,\circ\,(\,\D\otimes\mathrm{id}\,-\,2\<\nabla,\nabla>
  \,+\,\mathrm{id}\otimes\D\,)\,)\,(\,\psi\otimes\tilde\psi\,)
 $$
 which features three commuting operators $\D\otimes\mathrm{id}$,
 $\mathrm{id}\otimes\D$ and $\<\nabla,\nabla>$. Hence we are free to put
 these operators in arbitrary order upon exponentiation. For the ensuing
 calculations we need to choose a trivialization $\Phi$ of $\E$, which we
 think of as a family of homomorphisms $\Phi(x):\,\E_x\longrightarrow E$
 with $\Phi(0)=\mathrm{id}_E$, so that $\Phi\psi$ is a smooth section of
 $V\times E$ for every section $\psi$ of $\E$. With this in mind we find 
 \begin{eqnarray*}
  \left[\,e^{-z\D}\,(\,a(z)\,\psi\,)\,\right](0)
  &=&
  (\,\ev\,\circ\,e^{-z\D}\,\circ\,m\,)\,(\,a(z)\Phi^{-1}\otimes\Phi\psi\,)\\
  &=&
  \left(\,(\ev\otimes\ev)\,\circ\,e^{2z\<\nabla,\nabla>}\,\right)\,
  (\,e^{-z\D}(a(z)\Phi^{-1})\otimes e^{-z\D}(\Phi\psi)\,)\\
  &=&
  \<{(2z)^N\,(e^{-z\D}\,(a(z)\Phi^{-1})},{(e^{-z\D}\,(\Phi\psi))}>_{\S\,V^*}\\
  &=&
  \<{e^{z|x|^2}\,(2z)^N\,(e^{-z\D}\,(a(z)\Phi^{-1})},\Phi\psi>_{\S\,V^*}
 \end{eqnarray*}
 while $\left[e^{-z\QQ}\psi\right](0)\,=\,\<(\ev e^{-z\Phi\QQ\Phi^{-1}})^\#,
 \Phi\psi>_{\S\,V^*}$ by definition. As the scalar product $\<,>$ is
 non--degenerate on $\S\,V^*$ and $\psi$ can be chosen arbitrarily we
 conclude:
 \begin{equation}\label{link}
  (\,\ev e^{-z\Phi\QQ\Phi^{-1}}\,)^\#
  \quad=\quad e^{z|x|^2}\,(2z)^N\,(e^{-z\D}\,(a(z)\Phi^{-1}))\,.
 \end{equation}
 Evidently the operators $e^{z|x|^2}$ and $e^{-z\D}$ are invertible
 with inverses $e^{-z|x|^2}$ and $e^{z\D}$ respectively, thus we can
 solve equation (\ref{link}) uniquely for the infinite order jet of
 $a(z)$ leading to:
 \begin{eqnarray}\label{fi}
  \j^\infty_0(\,a(z)\Phi^{-1}\,)
  &=& e^{z\D}\,(2z)^{-N}\,e^{-z|x|^2}\,(\,\ev\,e^{-z\Phi\QQ\Phi^{-1}}\,)^\#
  \nonumber\\
  &=& e^{z\D}\,(2z)^{-N}\,(\ev\,e^{-z\Phi\QQ\Phi^{-1}}\,e^{z\D}\,)^\#\,.
 \end{eqnarray}
 Concerning the last step we remark that the polynomial $e^{-z|x|^2}$
 corresponds to the differential operator $e^{z\D}$ with constant
 coefficients. Hence it is possible to replace $e^{-z|x|^2}\,(\,
 \ev\,e^{-z\Phi\QQ\Phi^{-1}}\,)^\#$ by $(\,\ev\,e^{-z\Phi\QQ\Phi^{-1}}
 e^{z\D}\,)^\#$, but not by $(\,\ev\,e^{z\D}e^{-z\Phi\QQ\Phi^{-1}}\,)^\#$,
 because the latter would involve partial derivatives of the coefficients
 of $e^{-z\Phi\QQ\Phi^{-1}}$ in the origin. Nevertheless we cheated a
 little bit, because while certainly injective the operator $(2z)^N$
 is definitely not surjective. It is quite remarkable in itself that
 the value of the differential operator $e^{-z\Phi\QQ\Phi^{-1}}\,e^{z\D}$
 in the origin lies in the image of $(2z)^N$ by equation (\ref{link}),
 because its coefficients for the different powers $z^r,\,r\geq0,$ of $z$
 must have order less than or equal to $r$ to have this true! Hence the
 inverse $(2z)^{-N}$ is well defined and refers to the unique preimage
 of $(\ev\,e^{-z\Phi\QQ\Phi^{-1}}e^{z\D})^\#$ under $(2z)^N$.

 \begin{Definition}[The Difference Operator of a Generalized Laplacian]
 \label{dld}\endsummary
  Consider a generalized Laplacian $\Q$ and the flat Laplacian $\D$ acting
  on sections of $\E$ and the trivial vector bundle $V\times E$ over $V$
  respectively. Associated to $\Q$ and some smooth trivialization $\Phi$
  of $\E$ is the formal power series
  $$
   \DL(z)\quad:=\quad e^{-z(\Phi\QQ\Phi^{-1})}\,e^{z\D}
   \quad=\quad 1\,+\,z\,\DL_1\,+\,z^2\,\DL_2\,+\,z^3\,\DL_3\,+\,\cdots
  $$
  of differential operators on $V\times E$ called the difference operator.
  The values $\ev\DL_r$ of its coefficients have order less than or equal
  to $r$, so that $\ev\DL(z)$ lies in the image of $(2z)^N$.
 \end{Definition}

 At least in principle the recursion relation $\frac{d}{dz}\DL(z)\,=\,\DL(z)
 \D-(\Phi\Q\Phi^{-1})\DL(z)$ or equivalently $r\DL_r\,=\,\DL_{r-1}\D\,-\,
 (\Phi\Q\Phi^{-1})\DL_{r-1}$ for $r\geq 1$ together with $\DL_0=1$ allows
 us to calculate $\DL(z)$ to arbitrary order in $z$. However the composition
 $(\Phi\Q\Phi^{-1})\DL_{r-1}$ requires knowledge of the full operator
 $\DL_{r-1}$ and not only of its value at the origin. Note that the recursion
 relation in itself does not imply that the values of the $\DL_r,\,r\geq 0,$
 have order less than or equal to $r$, because the origin will be a singular
 point for the full operator $\DL(z)$ in general. Reformulating equation
 (\ref{fi}) according to the definition of $\DL(z)$ we arrive at the following
 theorem, which in a sense is the converse of Theorem \ref{iprop}:

 \begin{Theorem}[Characterization of the Heat Kernel Coefficients]
 \label{dlf}\endsummary
  The infinite order jet of the generating series $a(z)$ for the heat
  kernel coefficients of a generalized Laplacian $\Q$ acting on sections
  of a vector bundle $\E$ over $V$ is uniquely characterized by its
  intertwining property $\left[e^{-z\QQ}\psi\right](0)\,=\,\left[e^{-z\D}
  (a(z)\,\psi)\right](0)$. In fact the infinite order jet of any solution
  $a(z)$ to this equation is given in terms of some smooth trivialization
  $\Phi$ of $\E$ and the corresponding difference operator $\DL(z)$ by:
  $$
   \j^\infty_0(\,a(z)\Phi^{-1}\,)
   \quad=\quad e^{z\D}\;(2z)^{-N}\,(\ev\DL(z))^\#\,.
  $$
 \end{Theorem}

 \pskip
 The essence of Theorem \ref{dlf} is not particularly evident in the
 abstract formulation above and so we want to close this section expanding
 this formula for the jets of the heat kernel coefficients into powers of
 $z$. This exercise will show in particular that it is just as difficult
 to calculate the infinite order jets of the heat kernel coefficients $a_k,
 \,k\geq 0,$ as it is to calculate the value of the powers $\Q^r,\,r\geq 0,$
 of $\Q$ in $0$ using normal coordinates and the chosen trivialization
 $\Phi$ of $\E$. In other words calculating the powers $\Q^r$ explicitly
 is the problem of calculating the coefficients of the heat kernel
 expansion en nuce! Let us start by expanding the value $\ev\DL(z)$
 into powers of $z$ and homogeneous pieces
 $$
  (\,\ev\DL(z)\,)^\#\quad=\quad\sum_{r\geq s\geq 0}\,z^r\,\DL_{r,s}\,,
 $$
 where each $\DL_{r,s}\in\S^sV^*\otimes\mathrm{End}\,E$ is homogeneous
 polynomial of order $s$, and employ Theorem \ref{dlf} to conclude
 \begin{equation}\label{dlff}
  \j^\infty_0a(z)
  \quad=\quad\sum_{r\geq s\geq 2t\geq 0}
   \frac{z^t}{t!}\,\D^t\,(\,\frac{z^{r-s}}{2^s}\,\DL_{r,s}^\#\,)
  \quad=\quad\sum_{k,l\geq 0}\,z^k\,\sum_{t=0}^k\,\frac1{2^{l+2t}\,t!}\,
  \D^t\DL_{k+l+t,l+2t}^\#
 \end{equation}
 by reordering the sum and setting $k:=r-s+t$ and $l:=s-2t$. In particular
 the infinite order jet of the first coefficient $a_0$ picks up the principal
 symbols $\ev\DL_{l,l},\,l\geq 0,$ of the operators $\DL_l,\,l\geq 0,$ in
 the origin. Now by Definition \ref{dld} of the difference operator $\DL(z)$
 we find for the homogeneous pieces of its value
 $$
  \DL^\#_{r,s}\quad=\quad\sum_{\mu=0}^r\,(-1)^r\,\frac1{\mu!}|x|^{2\mu}\,
  \frac1{(r-\mu)!}(\ev\,(\Phi\Q\Phi^{-1})^{r-\mu})^\#_{s-2\mu}
 $$
 where $(\ev\,(\Phi\Q\Phi^{-1})^{r-\mu})^\#_{s-2\mu}\,\in\,\S^{s-2\mu}V^*
 \otimes\End E$ is the homogeneous piece of order $s-2\mu$ of the polynomial
 corresponding to $\ev\,(\Phi\Q\Phi^{-1})^{r-\mu}$. Inserting this expression
 into equation (\ref{dlff}) above we get a closed formula for the
 homogeneous pieces of the infinite order jets of the heat kernel
 coefficients $a_k,\,k\geq 0$. Though explicit this formula looks rather
 akward and so we refrain from writing it down. In fact the main point is
 to understand the fantastic cancellation properties of the $\Q_r^\#,\,
 r\geq 0,$ which make the $\ev\DL_r,\,r\geq 0,$ operators of degree less
 than or equal to $r$ in the first place! Without proper understanding of
 these cancellations we won't be able to take advantage of either formula
 to calculate the heat kernel coefficients explicitly.
\end{document}